\documentclass[letterpaper, 10 pt, conference]{ieeeconf}  

\IEEEoverridecommandlockouts
\usepackage[T1]{fontenc}
\usepackage[utf8]{inputenc}   

\usepackage{cite}
\usepackage{amsmath,amssymb,amsfonts,mathtools,bm,mathrsfs,soul,lipsum}
\usepackage{nicefrac}
\usepackage{algorithmic}
\usepackage{graphicx}
\usepackage{textcomp}
\usepackage{xcolor}
\usepackage[capitalise]{cleveref}
\usepackage{csquotes}

\newcommand{\tran}{^{\mathstrut\scriptscriptstyle\top}}
\DeclareMathOperator{\diag}{diag}

\DeclareMathOperator{\prob}{\mathsf{P}}

\newcommand{\numberOf}[1]{\ensuremath{N_{\mathrm{#1}}}}

\let\min\relax
\let\max\relax
\DeclareMathOperator*{\min}{{min}}
\DeclareMathOperator*{\max}{{max}}

\newcommand{\Max}{{\mathrm{max}}}
\newcommand{\Min}{{\mathrm{min}}}
\DeclareMathOperator*{\Minimize}{\mathbf{Minimize}}

\newcommand{\predHor}{N}

\DeclareMathOperator{\nodes}{\mathbf{nodes}}
\DeclareMathOperator{\stage}{\mathbf{stage}}
\DeclareMathOperator{\child}{\mathbf{child}}
\DeclareMathOperator{\ancestor}{\mathbf{anc}}

\newcommand{\One}{1}
\newcommand{\Zero}{0}
\newcommand{\Eye}{I}

\newcommand{\VaR}{\mathrm{V@R}}
\newcommand{\AVaR}{\mathrm{AV@R}}
\newcommand{\EVaR}{\mathrm{EV@R}}

\newcommand{\R}{\mathbb{R}}				
\newcommand{\N}{\mathbb{N}}				
\newcommand{\E}{\mathbf{E}}

\makeatletter
\newcommand{\subalign}[1]{%
  \vcenter{%
    \Let@ \restore@math@cr \default@tag
    \baselineskip\fontdimen10 \scriptfont\tw@
    \advance\baselineskip\fontdimen12 \scriptfont\tw@
    \lineskip\thr@@\fontdimen8 \scriptfont\thr@@
    \lineskiplimit\lineskip
    \ialign{\hfil$\m@th\scriptstyle##$&$\m@th\scriptstyle{}##$\crcr
      #1\crcr
    }%
  }
}
\makeatother

\definecolor{vgRed}{RGB}{193, 48, 24}
\definecolor{vgOrange}{RGB}{243, 111, 19}
\definecolor{vgYellow}{RGB}{235, 203, 56}
\definecolor{vgGreen}{RGB}{162, 185, 105}
\definecolor{vgLightBlue}{RGB}{13, 149, 188}
\definecolor{vgDarkBlue}{RGB}{6, 56, 81}

\usepackage{units}

\usepackage{tabularx}
\usepackage{csvsimple}
\usepackage{booktabs}
\usepackage{multirow}
\usepackage{multicol}
\usepackage{dcolumn}
\newcolumntype{d}[1]{D{.}{.}{#1}}

\newcommand{\noDecimal}[1]{\multicolumn{1}{c}{\centering #1}}

\usepackage{acronym}
\acrodef{res}[RES]{renewable energy sources}
\acrodef{pv}[PV]{photovoltaic}
\acrodef{mg}[MG]{microgrid}
\acrodef{ac}[AC]{alternating current}\acused{ac}
\acrodef{dc}[DC]{direct current}\acused{dc}
\acrodef{miqp}[MIQP]{mixed-integer quadratic program}
\acrodef{miqcp}[MIQCP]{mixed-integer quadratically-constrained program}
\acrodef{mpc}[MPC]{model predictive control}
\acrodef{arima}[ARIMA]{autoregressive integrated moving average}

\author{Jie Lei$^{1}$, Christian A. Hans$^{2}$ and Pantelis Sopasakis$^{3}$
\thanks{$^{1}$Jie Lei is with State Key Laboratory of Power Transmission Equipment \& System Security and New Technology at Chongqing University, China,
        {\tt\small jlei02@qub.ac.uk}.}%
\thanks{$^{2}$Christian A. Hans is with the Control Systems Group at Technische Universit\"{a}t Berlin, Germany,
        {\tt\small hans@control.tu-berlin.de}.}%
\thanks{$^{3}$Pantelis Sopasakis is with the School of EEECS and the i-AMS Centre at Queen's University Belfast, UK,
        {\tt\small p.sopasakis@qub.ac.uk}.}%
}

\title{\LARGE \bf Optimal operation of microgrids with risk-constrained state of charge%
}

\begin{document}

\maketitle

\begin{abstract}
        In this paper we present a stochastic scenario-based model predictive control (MPC) approach for the operation of islanded microgrids with high share of renewable energy sources.
        We require that the stored energy remains within given bounds with a certain probability using risk-based constraints as convex approximations of chance constraints.
        We show that risk constraints can generally be cast as conic constraints and, unlike chance constraints, can control both the number and average magnitude of constraint violations.
        Lastly, we demonstrate the risk-constrained stochastic MPC in a numerical case study.
\end{abstract}


\section{Introduction}

The high share of renewable energy sources in islanded \acp{mg} exposes
them to increased uncertainty due to the weather-dependent nature of
energy production. The pressing need to minimize the use of conventional
generators and maximize infeed from renewable sources while respecting
constraints on stored energy and power has led to a wide
adoption of \ac{mpc}
\cite{PBL+2013,HSB+2015,HSR+2020}.
\ac{mpc} additionally allows to make use of feed-forward information
using forecasting models (see, e.g., \cite{BJR2013}) of wind speed, solar irradiance and load demand.

Early works on deterministic \ac{mpc} \cite{TH2011,PRG2016,PBL+2013}
are giving way to methods that take the associated
uncertainty into consideration. Worst-case approaches \cite{HNRR2015,ZGG2013} can prove overly
conservative, especially in presence of higly uncertain wind, irradiance and
load.
Expectation-based (risk-neutral) stochastic \ac{mpc} formulations have been proposed
involving either (typically independent) processes with continuous distributions
\cite{GMV2015} or scenario-based formulations \cite{HSB+2015}.
Recently, a multistage risk-averse approach was proposed in \cite{HSR+2020} to
account for the fact that distributions are never known~exactly.

In stochastic \ac{mpc} formulations, constraints become stochastic too.
These can be imposed for all realizations of uncertainty \cite{HSR+2020,HSB+2015}.
A more appropriate and less conservative approach
is to require that the probability of constraints violations is sufficiently small.
Such \textit{chance-constrained} formulations are popular, e.g., in optimal
power flow problems~\cite{Xu2020}.
%
Chance constraints are generally nonconvex.
In certain cases, e.g., for certain continuous distributions and linear systems,
the inverse cumulative is known and chance constraints can be simplified \cite{GMV2015,Oldewurtel2008}.
The well-know framework of Nemirovski and Shapiro \cite{Nemirovski2007} for convex approximations of
chance constraints is widely known.
Some approximations have also been proposed such as
the stochastic tubes approach of \cite{SantosBonzanini2019} for constraint tightening,
and machine learning \cite{KLG2018}.

In this paper we propose a stochastic scenario-based \ac{mpc} scheme for the operation of islanded microgrids with a high share
of renewable sources.
We formulate a stochastic \ac{mpc} problem with chance constraints on the state of charge of storage units.
Moreover, we employ the approach of \cite{SSP2019a} to overapproximate chance constraints by convex conic constraints using coherent risk measures.
Lastly, we demonstrate the proposed chance-constrained MPC via realistic simulations with high renewable share using real-world irradiance and load data and time series forecast models.

\paragraph*{Notation}
Hereafter we denote by $\N$, $\N_0$ and $\R$ the sets of natural numbers, nonnegative integer
and real numbers respectively.
We denote the set of integers between $k$ and $k'$ by $\N_{[k, k']}$.
Let $\R_{> 0} = \{x\in\R{}\mid{} x>0\}$. Likewise define $\R_{\geq 0}$,
$\R_{\leq 0}$ and $\R_{< 0}$. The Euclidean norm is $\|{}\cdot{}\|_2$.
Given a set of indices $\N_{[k, k']}$, $[x_i]_{i\in\N_{[k, k']}}$ is shorthand for $[x_k, \cdots, x_{k'}]$.


\section{Microgrid model}
\label{sec:microgrid-model}

In what follows, some basics on scenario trees are introduced.
Moreover, the relation of model variables is discussed and a detailed overview over the \ac{mg} model is provided.

\subsection{Introduction on scenario trees}

In the context of this work, the probability distribution of load and available renewable infeed is assumed to be given in the form of a forecast scenario tree.

\begin{figure}
  \includegraphics{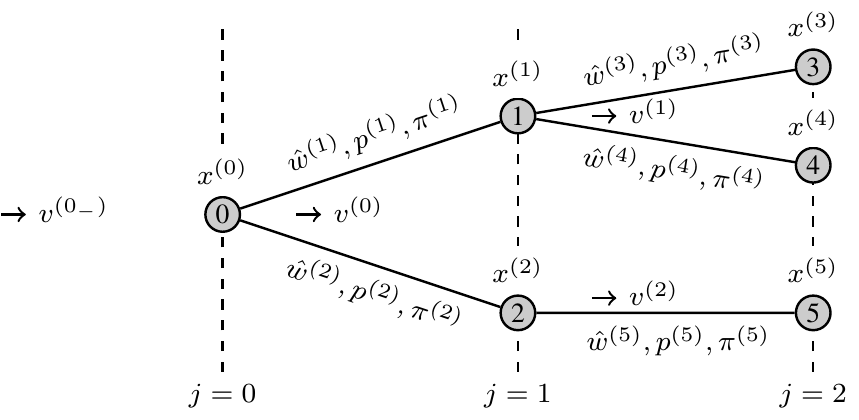}%
  \caption{Simple scenario tree from~\cite{HBR+2019}.
    In this tree,
    ${\nodes(0) = \{0\}}$,
    ${\nodes(1) = \{1, 2\}}$ and
    ${\nodes(2) = \{3, 4, 5\}}$.
    Moreover, the child nodes are
    $\child(0) = \{1, 2\}$,
    $\child(1) = \{3, 4\}$,
    ${\child(2) = \{5\}}$
    and consequently
    $\ancestor(3) = \ancestor(4) = 1$,
    $\ancestor(5) = 2$,
    $\ancestor(1) = \ancestor(2) = 0$.
  }%
  \label{fig:probabilityTree}
\end{figure}

A scenario tree is formed of $\numberOf{n}\in\N$ nodes.
Each node $i\in\N_{[0, \numberOf{n}-1]}$ is associated with prediction step $j\in\N_{[0, \predHor]}$, also referred to as stage, denoted as $\stage(i)$.
All nodes associated with stage $j$ are collected in the set $\nodes(j)$.
The number of nodes at stage $j$ is $n_j = |\nodes(j)|$.
Node $i = 0$ at stage $j=0$ is referred to as the root node, and the nodes at stage $\predHor$, i.e., the elements of $\nodes(\predHor)$, are referred to as leaf nodes.
Each node $i\in\nodes(j)$ at stage $j\in\N_{[1, \predHor]}$ is associated with an ancestor node $i_-\in\nodes(j-1)$ at stage $j-1$ which can be accessed via $i_- = \ancestor(i)$.
All nodes at stage $j$ that share a common ancestor $i\in\nodes(j-1)$ are referred to as child nodes of $i$ and are collected in $\child(i)\subseteq\nodes(j)$.
Each node $i$ is associated with a probability $\pi^{(i)}\in(0, 1]$.
For all stages $j\in\N_{[0,\predHor]}$, it holds that $\sum_{i\in\nodes(j)}\pi^{(i)} = 1$.
The probability of node $i\in\N_{[0, \numberOf{n}-1]}\setminus\nodes(\predHor)$ is linked to the probabilities of its children via $\pi^{(i)} = \sum_{i_+\in\child(i)} \pi^{(i_+)}$.
An example of a scenario tree is shown in \cref{fig:probabilityTree}.

\subsection{Generic \acl{mg} model}
\label{sec:modelling:generic}

\begin{table}[t]
  \centering
  \caption{Model-specific variables}
  \label{tab:modelSpecificVariables}
  \begin{tabular}{clcc}
    \toprule
    Symbol                & Explanation                           & Unit           & Size           \\
    \midrule
    $x$                   & Energy of storage units (state)       & $\unit{pu\,h}$ & $\numberOf{s}$ \\
    \midrule
    $u_{\mathrm{t}}$      & Control inputs of conventional units  & $\unit{pu}$    & $\numberOf{t}$ \\
    $u_{\mathrm{s}}$      & Control inputs of storage units       & $\unit{pu}$    & $\numberOf{s}$ \\
    $u_{\mathrm{r}}$      & Control inputs of renewable units     & $\unit{pu}$    & $\numberOf{r}$ \\
    $u$                   & Control inputs of all units           & $\unit{pu}$    & $\numberOf{u}$ \\
    $\delta_{\mathrm{t}}$ & Boolean control inputs of conv. units & ---            & $\numberOf{t}$ \\
    $v$                   & Vector of all control inputs          & ---            & $\numberOf{v}$ \\
    \midrule
    $w_{\mathrm{r}}$      & Uncertain available renewable power   & $\unit{pu}$    & $\numberOf{r}$ \\
    $w_{\mathrm{d}}$      & Uncertain load                        & $\unit{pu}$    & $\numberOf{d}$ \\
    $w$                   & Vector of all uncertain inputs        & $\unit{pu}$    & $\numberOf{w}$ \\
    \midrule
    $p_{\mathrm{t}}$      & Power of conventional units           & $\unit{pu}$    & $\numberOf{t}$ \\
    $p_{\mathrm{s}}$      & Power of storage units                & $\unit{pu}$    & $\numberOf{s}$ \\
    $p_{\mathrm{r}}$      & Power of renewable units              & $\unit{pu}$    & $\numberOf{r}$ \\
    $p$                   & Power of all units                    & $\unit{pu}$    & $\numberOf{u}$ \\
    $p_{\mathrm{e}}$      & Power over transmission lines         & $\unit{pu}$    & $\numberOf{e}$ \\
    \bottomrule
  \end{tabular}
\end{table}

We consider an \ac{mg} that comprises an arbitrary finite number of conventional and renewable generators, storage units and loads.
These components are connected to each other via \ac{ac} transmission lines.
With the variables summarized in \cref{tab:modelSpecificVariables}, the behavior of the \ac{mg} is modeled for all nodes $i_+\in\N_{[1, \numberOf{n}-1]}$ with $i=\ancestor(i_+)$ by
\begin{subequations}\label{eq:genericModel}
  \begin{align}
    x^{(i_+)} & = A x^{(i)} + B p^{(i_+)}, \label{eq:genericModel:dynamics}                               \\
    h_x       & \leq H_x x^{(i_+)}, \label{eq:genericModel:stateLimit}                                    \\
    p^{(i_+)} & = f_p(v^{(i)}, w^{(i_+)}), \label{eq:genericModel:power}                                  \\
    h_{vp}    & \leq H_{vp} \big[\begin{matrix} {v^{(i)}}\tran & {p^{(i_+)}}\tran \end{matrix}\big]\tran. \label{eq:genericModel:powerLimits}
  \end{align}
\end{subequations}
Here, $x^{(i)}$ represents the state, i.e., the energy stored in storage units at node $i$.
Moreover,
$u^{(i)} = [\begin{matrix} {u_{\mathrm{t}}^{(i)}}\tran & {u_{\mathrm{s}}^{(i)}}\tran & {u_{\mathrm{r}}^{(i)}}\tran \end{matrix}]\tran$
is a vector of real-valued control inputs and
$\delta_{\mathrm{t}}^{(i)}\in\{0, 1\}^{\numberOf{t}}$
a vector of Boolean inputs that indicates whether conventional unit $l\in\N_{[1, \numberOf{t}]}$ is enabled ($\delta_{\mathrm{t},l}^{(i)} = 1$) or disabled ($\delta_{\mathrm{t},l}^{(i)} = 0$).
The control inputs are collected in
$v^{(i)} = [\begin{matrix} {u^{(i)}}\tran &  {\delta^{(i)}}\tran \end{matrix}]\tran$.
Finally, the uncertain input at node ${i_+\in\child(i)}$ is
$w^{(i_+)} = [\begin{matrix} {w_{\mathrm{r}}^{(i_+)}}\tran & {w_{\mathrm{d}}^{(i_+)}}\tran \end{matrix}]\tran$
and the units' power is
${p^{(i_+)} = [\begin{matrix} {p_{\mathrm{t}}^{(i_+)}}\tran & {p_{\mathrm{s}}^{(i_+)}}\tran & {p_{\mathrm{r}}^{(i_+)}}\tran \end{matrix}]\tran}$.

\begin{figure}
  \centering
  \includegraphics{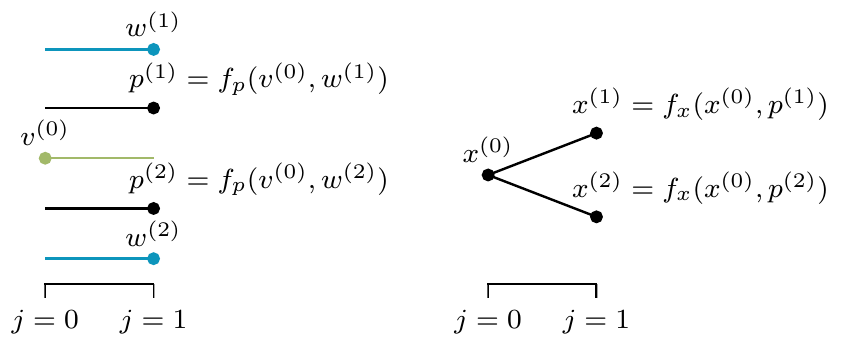}
  \caption{Relation of variables in scenario tree. Source:~\cite{HBR+2019}.}
  \label{fig:probabilitytreeVariables}
\end{figure}

The control input $v^{(i)}$ at node $i\in\N_{[0, \numberOf{n}]}\setminus\nodes(\predHor)$ is present between time instants $j=\stage(i)$ and $j+1$.
The uncertain input $w^{(i_+)}$ and the power $p^{(i_+)}$ at nodes $i_+\in\child(i)$ are associated with the same time interval.
The reason for this is that power and uncertain input are associated with the state $x^{(i_+)}$ that they result in.

An example of this relation is illustrated in \cref{fig:probabilitytreeVariables}.
Here, the control input between instants $j=0$ and $j=1$ is $v^{(0)}$.
In the example, two realizations of the uncertain input are predicted, $w^{(1)}$ and $w^{(2)}$.
Thus, for the same control input $v^{(0)}$, different values of load and available renewable infeed can occur.
The power of the units, which is collected in $p^{(1)}$ and $p^{(2)}$, changes with the control input and the uncertain input.
Similarly, the state changes with these power values.
As shown in \cref{fig:probabilitytreeVariables}, $x^{(1)}$ is a function of $x^{(0)}$ and $p^{(1)}$ and $x^{(2)}$ a function of $x^{(0)}$ and $p^{(2)}$.

Note that by introducing additional free variables, \eqref{eq:genericModel:power} can be transformed into a set of affine equality and inequality constraints (see also \cref{sec:energyRelatedLimits} and \cite{HSR+2020}).
Thus, \eqref{eq:genericModel} only comprises affine constraints and can therefore be used to formulate \aclp{miqp}.

In what follows, we will discuss the equations used to formulate the generic model \eqref{eq:genericModel}.
We start with the dynamics.

\subsection{Dynamics}
The dynamics of the stored energy $x^{(i_+)}$ are modeled using a discrete-time state model.
With
${A = \Eye_{\numberOf{s}}}$
and
${B = [\begin{matrix} \Zero_{\numberOf{s}\times\numberOf{t}} & -T_{\mathrm{s}} \Eye_{\numberOf{s}} & \Zero_{\numberOf{s}\times\numberOf{r}} \end{matrix}]}$
the model reads
\begin{equation} \label{eq:microgridModel:energy:dynamics}
  x^{(i_+)} = A x^{(i)} + B p^{(i_+)} \quad \text{with} \quad x(0) = x_0
\end{equation}
and $i= \ancestor(i_+)$. This is precisely \Cref{eq:genericModel:dynamics}.

\subsection{Energy-related limits}
\label{sec:energyRelatedLimits}
The stored energy at all nodes $i_+\in\N_{[1, \numberOf{n}]}$ is limited by
\begin{equation}
  x^{\min} \leq x^{(i_+)} \leq x^{\max}
\end{equation}%
with $x^{\min}\in\R^{\numberOf{s}}_{\geq 0}$ and $x^{\max}\in\R^{\numberOf{s}}_{\geq 0}$.
From these inequalities, \eqref{eq:genericModel:stateLimit} can be easily deduced.

\subsection{Power of units}

The units' power is not only affected by the power setpoints but also by the uncertain load and available renewable infeed.
This is taken into account by considering different realisations of the uncertain input $w^{(i_+)}$ for each control input $v^{(i)}$.
The effects of the control input and the uncertain input on the power enter the generic model \eqref{eq:genericModel} via function $f_p$ which is implicitly defined by constraints \eqref{eq:microgrid:power:equilibrium}--\eqref{eq:microgrid:power:gridForming}.

\subsubsection{Transmission network}
A power equilibrium of generation, consumption and storage power must be ensured at all times.
This can be modelled by the equality constraint
\begin{equation}\label{eq:microgrid:power:equilibrium}
  \One_{\numberOf{t}}\tran  p_{\mathrm{t}}^{(i_+)} + \One_{\numberOf{s}}\tran  p_{\mathrm{s}}^{(i_+)} + \One_{\numberOf{r}}\tran  p_{\mathrm{r}}^{(i_+)} + \One_{\numberOf{d}}\tran  w_{\mathrm{d}}^{(i_+)} = 0.
\end{equation}

\subsubsection{Renewable units}
The power of the renewable units $p_{\mathrm{r}}^{(i_+)}$ cannot exceed the weather-dependent available renewable infeed $w_{\mathrm{r}}^{(i_+)}$.
The control input $u_{\mathrm{r}}^{(i)}$ allows to limit $p_{\mathrm{r}}^{(i_+)}$ below $w_{\mathrm{r}}^{(i_+)}$.
If the available renewable power $w_{\mathrm{r},l}^{(i_+)}$ of unit $l\in\N_{[1,\numberOf{r}]}$ is below the power setpoint $u_{\mathrm{r},l}^{(i)}$, then the infeed of this unit equals the available renewable power.
If the available renewable power $w_{\mathrm{r},l}^{(i_+)}$ is above the power setpoint $u_{\mathrm{r},l}^{(i)}$, then the infeed of this unit equals the power setpoint.
This can be modelled via
\begin{equation} \label{eq:microgridModel:power:renewable}
  p_{\mathrm{r}}^{(i_+)} = \min (u_{\mathrm{r}}^{(i)}, w_{\mathrm{r}}^{(i_+)}).
\end{equation}

\subsubsection{Power sharing of grid-forming units}
The lower control layers of the \ac{mg} are assumed to be designed such that fluctuations of load and renewable infeed are distributed among grid-forming storage and conventional units in a proportional manner.
This ensures that \eqref{eq:microgrid:power:equilibrium} holds in presence of uncertain load and renewable infeed.
It can be implemented using, for example, decentralized droop control \cite{SHK+2017}.
Typically, power sharing is implemented on a much faster time scale than operation control.
Nevertheless, power sharing needs to be considered in the constraints as it links the power of the conventional and storage units to the fluctuations of the uncertain input.
Using the additional variable $\mu^{(i_+)}\in\R$, it can be included in the \ac{mpc} formulation by \cite{HSR+2020}
\begin{subequations}\label{eq:microgrid:power:gridForming}
  \begin{align}
    {K}_\mathrm{s} (p_\mathrm{s}^{(i_+)} - u_\mathrm{s}^{(i)}) & = \mu^{(i_+)},                         \\
    {K}_\mathrm{t} (p_\mathrm{t}^{(i_+)} - u_\mathrm{t}^{(i)}) & = \mu^{(i_+)} \delta_\mathrm{t}^{(i)}.
  \end{align}
\end{subequations}

As indicated in \cite{HSR+2020}, \eqref{eq:microgrid:power:equilibrium}--\eqref{eq:microgrid:power:gridForming} can be equivalently expressed by affine equality and inequality constraints using additional decision variables and the so-called Big-M method (see, e.g.,~\cite{BM1999}).
Consequently, $f_p(v^{(i)}, w^{(i_+)})$ in \eqref{eq:genericModel:power} can be used to formulate mixed-integer optimization problems.

\subsection{Power-related limits}

The limits on power and power setpoints can be divided into limits on units and limits on transmission lines.
In detail, \eqref{eq:genericModel:powerLimits} is composed of the following inequalities.

\subsubsection{Renewable units}
Power and control input are limited by $p_{\mathrm{r}}^{\min}\in\R^{\numberOf{r}}_{\geq 0}$ and $p_{\mathrm{r}}^{\max}\in\R^{\numberOf{r}}_{\geq 0}$, i.e.,
\begin{subequations}\label{eq:microgridModel:renewable:limits}
  \begin{alignat}{2}
    p_{\mathrm{r}}^{\min} & \leq p_{\mathrm{r}}^{(i_+)} &  & \leq p_{\mathrm{r}}^{\max}, \\
    p_{\mathrm{r}}^{\min} & \leq ~ u_{\mathrm{r}}^{(i)} &  & \leq p_{\mathrm{r}}^{\max}.
  \end{alignat}
\end{subequations}

\subsubsection{Conventional units}
If a unit is disabled, then power and setpoint are zero,
otherwise power and setpoint are limited by a minimum and a maximum value.
\begin{subequations}\label{eq:microgridModel:conventional:limits}%
  With ${p_{\mathrm{t}}^{\min}\in\R^{\numberOf{t}}_{\geq 0}}$ and $p_{\mathrm{t}}^{\max}\in\R^{\numberOf{t}}_{\geq 0}$, this can be expressed by
  \begin{alignat}{2}
    \diag(p_\mathrm{t}^{\min}) \delta_{\mathrm{t}}^{(i)} & \leq p_{\mathrm{t}}^{(i_+)} &  & \leq \diag(p_\mathrm{t}^{\max}) \delta_{\mathrm{t}}^{(i)}, \\
    \diag(p_\mathrm{t}^{\min}) \delta_{\mathrm{t}}^{(i)} & \leq ~ u_{\mathrm{t}}^{(i)} &  & \leq \diag(p_\mathrm{t}^{\max}) \delta_{\mathrm{t}}^{(i)}.\end{alignat}%
\end{subequations}%

\subsubsection{Storage units}
The power and the power setpoints are limited by $p_{\mathrm{s}}^{\min}\in\R^{\numberOf{s}}_{\leq 0}$ and $p_{\mathrm{s}}^{\max}\in\R^{\numberOf{s}}_{\geq 0}$, i.e.,
\begin{subequations}\label{eq:microgridModel:storage:limits}
  \begin{alignat}{2}
    p_{\mathrm{s}}^{\min} & \leq p_{\mathrm{s}}^{(i_+)} &  & \leq p_{\mathrm{s}}^{\max}, \\
    p_{\mathrm{s}}^{\min} & \leq ~ u_{\mathrm{s}}^{(i)} &  & \leq p_{\mathrm{s}}^{\max}.
  \end{alignat}
\end{subequations}

\subsubsection{Transmission network}
The transmission network is included using the linear \ac{dc} power flow approximations for \ac{ac} grids (see, e.g. \cite{HBR+2019}).
As the power equilibrium is already implemented via \eqref{eq:microgrid:power:equilibrium}, we can directly deduce the power flowing over the transmission lines $p_{\mathrm{e}}^{(i_+)}$ via
\begin{equation*}
  p_{\mathrm{e}}^{(i_+)} = F\,
  \big[
    \begin{matrix}
      p_{\mathrm{t}}^{(i_+)} &
      p_{\mathrm{s}}^{(i_+)} &
      p_{\mathrm{r}}^{(i_+)} &
      w_{\mathrm{d}}^{(i_+)}
    \end{matrix}
    \big]\tran
\end{equation*}
with $F\in\R^{\numberOf{e}\times(\numberOf{u}+\numberOf{d})}$.
Naturally, $p_{\mathrm{e}}^{(i_+)}$ is bounded, i.e.,
\begin{equation}\label{eq:microgridModel:line:limits}
  p_{\mathrm{e}}^{\min} \leq F\,
  \big[
    \begin{matrix}
      p_{\mathrm{t}}^{(i_+)} &
      p_{\mathrm{s}}^{(i_+)} &
      p_{\mathrm{r}}^{(i_+)} &
      w_{\mathrm{d}}^{(i_+)}
    \end{matrix}
    \big]\tran \leq p_{\mathrm{e}}^{\max}
\end{equation}
with
$p_{\mathrm{e}}^{\min}\in\R_{\leq 0}^{\numberOf{e}}$
and
$p_{\mathrm{e}}^{\max}\in\R_{\geq 0}^{\numberOf{e}}$.

This completes the introduction of the control-oriented \ac{mg} model.
Based on this section, we will now formulate a cost function for islanded \ac{mg} with high renewable share.


\section{Operating costs}
\label{sec:operatingCosts}

The operating cost of the microgrid is motivated by \cite{HSB+2015}.
It is composed of
(i) fuels costs of the conventional units, $\ell_{\mathrm{t}}^{\mathrm{f}}$,
(ii) switching costs of the conventional units, $\ell_{\mathrm{t}}^{\mathrm{s}}$, and
(iii) costs incurred by limiting potential renewable infeed, $\ell_\mathrm{r}$.
For scenario trees of the form described in \Cref{sec:microgrid-model}, the cost associated with node $i_+\in\N_{[1,\numberOf{n}]}$ is
\begin{multline}
	\label{eq:operatingCosts:economic}
	\ell (v^{(i)}, v^{(i_-)}, p^{(i_+)}) = \big(
	\ell_{\mathrm{t}}^{\mathrm{f}} (v^{(i)}, p^{(i_+)}) +
	\\
	\ell_{\mathrm{t}}^{\mathrm{s}}(v^{(i)}, v^{(i_-)}) +
	\ell_\mathrm{r}(p^{(i_+)})\big)\gamma^{\stage(i_+)},
\end{multline}
with $i=\ancestor(i_+)$ and $i_-=\ancestor(i)$.
Here, discount factor $\gamma\in(0, 1]$ is used to put an emphasis on near future decisions.

In \eqref{eq:operatingCosts:economic}, decision variables associated with different nodes are used.
The reason for this is that the power $p^{(i_+)}$ depends on control input $v^{(i)}$, $i=\ancestor(i_+)$.
In addition, the switching costs associated with $v^{(i)}$ depend on the Boolean input $\delta_\mathrm{t}^{(i_-)}$ at the previous time instant which is part of $v^{(i_-)}$.

\begin{subequations}
	The fuel cost of conventional units is approximately \cite{JMK2012}
	\begin{equation}
		\ell_{\mathrm{t}}^{\mathrm{f}} (v^{(i)}, p^{(i_+)})
		= c_{\mathrm{t}}\tran \delta_{\mathrm{t}}^{(i)} + {c_{\mathrm{t}}^\prime}\tran p^{(i_+)}_{\mathrm{t}} + \|\diag(c_{\mathrm{t}}^{\prime\prime}) p^{(i_+)}_{t} \|_2^2
	\end{equation}
	with weights
	$c_{\mathrm{t}}\in\R^T_{>0}$,
	$c_{\mathrm{t}}^{\prime}\in\R^T_{>0}$, and
	$c_{\mathrm{t}}^{\prime\prime}\in\R^T_{>0}$.
	Moreover, the cost incurred by switching the conventional generators is
	\begin{equation}\label{eq:cost-function:switch}
		\ell_{\mathrm{t}}^{\mathrm{s}}(v^{(i)}, v^{(i_-)}) =  \|\diag(c_{\mathrm{t}}^{\mathrm{s}}) (\delta_\mathrm{t}^{(i_-)} - \delta_\mathrm{t}^{(i)} ) \|_2^2
	\end{equation}
	with weight
	$c_{\mathrm{t}}^{\mathrm{s}}\in\R^T_{>0}$.
	For nodes $i_+\in\N_{[1, \numberOf{n}-1]}$ with $\ancestor(i_+) = 0$, \eqref{eq:cost-function:switch} becomes $\ell_{\mathrm{t}}^{\mathrm{s}}(v^{(0)}, v^{(0_-)})$.
	Here, $v^{(0_-)}$ is the input applied at the previous execution of the controller.

	The goal in the operation of renewable units is to use as much weather-dependent available power as possible.
	Using $p_{\mathrm{r}}^{\Max}$ from \cref{sec:microgrid-model} and weight $c_{\mathrm{r}}\in\R^R_{>0}$, this can be encoded into the cost function via
	\begin{equation}
		\ell_\mathrm{r}(p^{(i_+)}) = \|\diag(c_{\mathrm{r}}) (p_{\mathrm{r}}^{\max} - p^{(i_+)}_{\mathrm{r}})\|_2^2.
	\end{equation}
\end{subequations}

\begin{subequations}\label{eq:costVectorAtStageJPlusOne}
	In favour of a more compact notation, we introduce the cost variable $Z^{(i_+)}$ for every node $i_+\in\N_{[1, \numberOf{n}-1]}$, i.e.,
	\(
	Z^{(i_+)} = \ell(v^{(i)}, v^{(i_-)}, p^{(i_+)}).
	\)
	The cost values at stage $j$ can be collected in
	\(
	Z_{j} = [Z^{(i_+)}]_{i_+ \in \nodes(j)},
	\)
	which defines a random variable on the probability space $\nodes(j)$.
\end{subequations}
Using the probabilities of stage $j$, $\pi_j = [\pi^{(i_+)}]_{i_+ \in \nodes(j)}$, the expected cost at stage $j$ is
\begin{equation}
	\E^{\pi_j}(Z_j) = \pi_j\tran Z_j.
\end{equation}
The total cost along the prediction horizon, $Z\in\R$, is the sum of expected costs over all stages of the tree, i.e.,
\begin{equation}\label{eq:expected_predicted_cost}
	Z = \textstyle\sum\limits_{j=1}^\predHor \E{\pi_j}(Z_j) = \textstyle\sum\limits_{i_+=1}^{\numberOf{n}-1} \pi^{(i_+)} Z^{(i_+)}.
\end{equation}
Using the cost and the \ac{mg} model from \cref{sec:microgrid-model}, we can now formulate a chance-constrained \ac{mpc} problem for islanded \acp{mg}.

\newcommand{\chance}{\alpha}

\section{Risk-constrained model predictive control}
In this section we shall state the chance-constrained
and risk-constrained \ac{mpc} problems.

\subsection{Chance-constrained optimal control formulations}\label{sec:chance_constraints}
Consider the following stochastic optimal control problem with horizon $\predHor$
\begin{align}
    \label{eq:optimal-control-problem}
    \Minimize_{\bm{v}, \bm{x}, \bm{p}}\,Z\tag{P}
\end{align}
subject to the system dynamics and constraints in Equations
\eqref{eq:genericModel} for all $i_+\in\N_{[1, \numberOf{n}-1]}$ with $i=\ancestor(i_+)$,
which consists in minimizing the expected
value of the cost defined in \eqref{eq:expected_predicted_cost}.
The minimization is carried out over $\bm{x}=[x^{(i)}]_{i\in\N_{[0, \numberOf{n}-1]}}$,
$\bm{v} = [v^{(i)}]_{i\in\N_{[0, \numberOf{n}-1]}\setminus\nodes(\predHor)}$
and $\bm{p} = [p^{(i)}]_{i\in\N_{[1, \numberOf{n}-1]}}$
and given that the state at the root node, $x^{(0)}$, is equal to the
measured state.
We need to impose that the state of charge of each storage unit $s\in\N_{[1, \numberOf{s}]}$
at each stage $j$, $x_{s,j}=[x_{s}^{(i)}]_{i\in\nodes(j)}$,
which is a random variable, remains bounded between \(\tilde{x}^{\Min} \geq x^\Min\)
and \(\tilde{x}^{\Max} \leq x^\Max\) in a probabilistic sense.
These tighter bounds serve to reduce the range of the state of charge which typically
comes with positive effect on the aging of storage units (see, e.g., \cite{AZKC2016}).

To require that $x_{s,j}$ satisfies the constraints for all possible
realizations of the uncertain renewable infeed and load can be overly conservative.
A possible alternative is to require that the probability of constraint
satisfaction at every storage unit $s$ is adequately high, i.e.,
\begin{equation}\label{eq:chance_constraints}
    \prob_j\left[x_{s,j} \notin [\tilde{x}^{\Min}, \tilde{x}^{\Max}]\right] \leq \chance,
\end{equation}
for all stages $j\in\N_{[1,\predHor]}$, for some $\chance\in[0, 1]$,
where \(\prob_j\) is the probability measure of $\nodes(j)$ associated with probability vector $\pi_j$.
This condition can be rewritten using the value-at-risk operator
of a random variable $X$ on the probability space $\nodes(j)$,
\begin{equation}
    \VaR_{\chance}[X] \coloneqq \inf\{t {}:{} \prob_j[X > t] \leq \chance\}.
\end{equation}
In fact, the above probabilistic constraints are equivalent to
\begin{equation}\label{eq:var_constraints}
    \VaR_{\chance}\left[d_j(x_{s,j})\right] \leq 0,
\end{equation}
for all $j\in\N_{[1, \predHor]}$, where $d_j$ is the
distance-to-set function defined as
\(
d_j(x)
{}={}
\min_{y\in [\tilde{x}_{s}^{\Min}, \tilde{x}_{s}^{\Max}]}|x-y|.
\)
One-sided constraints of the form
\(
\prob_j\left[x_{s,j} > \tilde{x}^{\Max}\right] \leq \chance,
\)
can be imposed by considering functions of the form
\(
{d_j^{\Max}(x)
{}={}
        \min_{y \leq \tilde{x}_{s}^{\Max}}|x-y|.}
\)

Such probabilistic constraints are nonconvex,
call for additional binary variables
and can lead to computationally intractable optimization problems,
while they limit the frequency of violations, but not their magnitude.
Indeed, in order to impose the constraint of \eqref{eq:chance_constraints}
we  define the binary variable $\tau^{(i)}_s\in\{0, 1\}^{n_j}$
which is such that
\begin{equation}\label{eq:tau}
    \tau^{(i)}_s = 0
        {}\iff{}
    d_j(x_{s}^{(i)}) \leq  0,
\end{equation}
for all $i\in\nodes(j)$.
Using a standard big-M relaxation, we choose $M \geq \max(x^{\Max} - x^\Min)$,
$m\leq 0$ and a small tolerance $\epsilon>0$ to rewrite
\eqref{eq:tau} as
\begin{align}
    \epsilon + (m-\epsilon)\tau_j^{(i)} \leq d_j \leq M(1-\tau_j^{(i)}),
\end{align}
for all $i\in\nodes(j)$.
Then, for each storage unit $s$,
\begin{multline}
    \prob_j[d_j(x_{s}^{(i)}) > 0]\leq \chance
    \iff
    \hspace{-0.7em}
    \textstyle \sum\limits_{i\in\nodes(j)}\hspace{-0.9em}\pi^{(i)}\tau_{s}^{(i)} \geq 1-\chance,
\end{multline}
which, by defining $\tau_{s,j}=[\tau_{s}^{(i)}]_{i\in\nodes(j)}$,
can be equivalently written as
\(
\pi_j^\top \tau_{s,j} \geq 1-\chance.
\)
In summary, the above chance constraints are equivalent to
\begin{equation}\label{eq:chance-constraints-binary}
    \begin{aligned}
        \tau_{s,j}^{(i)} \in \{0, 1\}, \quad d_j(x_{s,j}^{(i)}) \leq \xi_{s}^{(i)}, \quad
        \pi_j^\top \tau_{s,j} \geq 1-\chance,
        \\
        \epsilon + (m-\epsilon)\tau_{s,j}^{(i)} \leq \xi_{s}^{(i)} \leq M(1-\tau_j^{(i)}),
    \end{aligned}
\end{equation}
for $i\in\nodes(j)$, where we have introduced a relaxation of
$d_j(x_{s,j}^{(i)})$ by introducing the
auxiliary variables $\xi_s^{(i)}$. Note that in order to impose chance constraints
we need to introduce as many binary variables as the nodes of the tree times the number
of storage units. Overall, the chance-constrained optimal control problem reads
\begin{align}
    \label{eq:chance-constrained-optimal-control-problem}
    \Minimize_{\bm{v}, \bm{x}, \bm{p}, \bm{\xi}, \bm{\tau}}\,Z\tag{P$_{\text{cc}}$}
\end{align}
subject to the system dynamics and constraints in
\eqref{eq:genericModel} and \eqref{eq:chance-constraints-binary}
for all $i_+\in\N_{[1, \numberOf{n}-1]}$ with $i=\ancestor(i_+)$.
Note that the minimization is taken over \(\bm{v}\), \(\bm{x}\), \(\bm{p}\)
and the auxiliary variables
\({\bm{\xi} = [[\xi^{(i)}_s]_{s\in\N_{[1, \numberOf{s}]}}]_{i\in\N_{[1, \numberOf{n}-1]}}}\),
\({\bm{\tau} = [[\tau^{(i)}_s]_{s\in\N_{[1, \numberOf{s}]}}]_{i\in\N_{[1, \numberOf{n}-1]}}}\).

An alternative approach is to use risk-based constraints, which
can control both the occurrence and the magnitude of violations.

\subsection{Risk measures}
A risk measure is an operator that maps random variables to characteristic
values. The risk of a random cost represents an ``equivalent'' fixed
cost value. The expectation and maximum operators are examples of risk measures.

A risk measure $\rho$ on a finite probability space $(\Omega, \prob)$,
with $\Omega = \{i\}_{i=1}^{n}$,
is considered to be well-behaving if it is
(i) convex, that is, for all random variables
$X_1$, $X_2$ and $\lambda \in [0,1]$, it is
\(
\rho(\lambda X_1 + (1-\lambda)X_2)
\leq
\lambda \rho(X_1) + (1-\lambda) \rho(X_2),
\)
(ii) monotone, in the sense that $\rho(X_1) \leq \rho(X_2)$ whenever
$\prob[X_1 > X_2] = 0$,
(iii) translation equi-variant, i.e.,
\(
\rho(X + c) = c + \rho(X)
\)
for all constants $c\in\R$, and
(iv) positive homogeneous, that is
\(
\rho(aX) = a \rho(X),
\)
for all $a\geq 0$.
Risk measures that satisfy these requirements are
called \textit{coherent} and can be represented as
\begin{equation}
    \label{eq:risk-measure-representation}
    \rho[X] = \sup_{\mu \in \mathcal{A}}\E^{\mu}[X],
\end{equation}
where $\E^{\mu}$ is the expectation operator with respect to a probability
$\mu$ and $\mathcal{A}$ is a convex set of probabilities
called the \textit{ambiguity set} of $\rho$ \cite[Theorem 6.5]{SDR2014}.

A widely used coherent risk measure is the
\textit{average value-at-risk at level $\chance$}
denoted by $\AVaR_{\chance}$;
its ambiguity set over a finite, $n$-dimensional,
probability space  with probability vector $\pi\in\R^n$,
is the polytope
\begin{equation}
    \mathcal{A}_\chance =
    \left\{
    \mu{}\in{}\R^n
    {}:{}
    \textstyle\sum_{i=1}^{n}\mu_i = 1,
    0\leq \chance\mu \leq \pi
    \right\}.
\end{equation}
Note that $\mathcal{A}_0$ coincides with the whole probability simplex
and $\AVaR_0[X] {}={} \max_{i=1,\ldots, n}\{X_i {}:{} \pi_i \neq 0\}$,
while $\mathcal{A}_1 = \{\pi\}$, therefore $\AVaR_1[X] = \E^{\prob}[X]$.

A noteworthy property of the average value-at-risk is that it is
a tight convex overapproximation of the
value-at-risk, that is $\AVaR_{\chance}[X] \geq \VaR_{\chance}[X]$,
therefore the probabilistic constraints of \eqref{eq:var_constraints}
are satisfied if
\begin{equation}\label{eq:avar-constraints}
    \AVaR_{\chance}\left[
        d_j(x_{s,j})
        \right] \leq 0,
\end{equation}
for all $j\in\N_{[1, \predHor]}$.
Unlike $\VaR$-based constraints, $\AVaR$-based constraints are convex
and do not lead to overly cumbersome optimization problems.
For given $\chance \in [0, 1]$, any $\VaR_\chance$-bounding risk measure
with $\rho[Z] \geq \VaR_{\chance}[Z]$ will
imply the satisfaction of the original probabilistic constraints.

The entropic value-at-risk, $\EVaR_\chance$, is another example of such a risk measure \cite{AhmadiJavid2011}; its ambiguity set over a finite, $n$-dimensional, probability space with probability vector $\pi$ is the set of probability vectors, $\mu\in\R^n$, such that $\mathrm{D}_{\rm KL}(\mu {}\Vert{} \pi) \leq -\ln \chance$, where \(\mathrm{D}_{\rm KL}\) is the Kullback-Leibler divergence.

Lastly, the level of risk aversion can be determined from available data using statistical methods \cite{SSP2019}.

\subsection{Risk constraints}\label{sec:risk_constraints}
Given a coherent $\VaR$-bounding risk measure \(\rho_j\)
defined on the probability space $\nodes(j)$, which is equipped with the
probability vector $\pi_j$, probabilistic constraints of the form given in
\eqref{eq:var_constraints} can be overapproximated the
convex constraints
\begin{equation}\label{eq:stagewise-risk-constraints}
    \rho_j[C_j] \leq 0,
\end{equation}
where
\(
C_j
{}={}
d_{j}(x_j),
\) for $j\in\N_{[1, \predHor]}$.
In this section we have dropped the index $s$ (cf. \Cref{eq:avar-constraints})
for the sake of simplicity.
The associated risk-constrained optimal control problem is
\begin{align}
    \label{eq:risk-constrained-optimal-control-problem}
    \Minimize_{\bm{v}, \bm{x}, \bm{p}}\,Z\tag{P$_{\text{rc}}$}
\end{align}
subject to \eqref{eq:stagewise-risk-constraints} for all $j\in\N_{[1, \predHor]}$
and the system dynamics and constraints in \eqref{eq:genericModel} for all
$i_+\in\N_{[1, \numberOf{n}-1]}$ with ${i=\ancestor(i_+)}$.

Since \(\rho_j\) is a coherent risk measure, there is
a closed convex set $\mathcal{A}_{\chance, j}$ such that
\begin{equation}
    \rho_j[C_j]
    {}={}
    \max_{\mu \in \mathcal{A}_{\chance, j}} \E^{\mu}[C_j]
    {}={}
    \max_{\mu \in \mathcal{A}_{\chance, j}} \mu\tran C_j.
\end{equation}
The ambiguity set $\mathcal{A}_{\chance, j}$ can generally be written in
a conic form as follows
\begin{equation}
    \mathcal{A}_{\chance, j}
    {}={}
    \left\{\mu
    \in \R^{n_j}
    \left|
    \begin{array}{l}
        \exists \nu \in \R^{r_j}{}\text{ such that }
        \\
        b_{\chance,j} {-} E_{\chance,j}\mu {-} F_{\chance,j}\nu {}\in{} \mathcal{K}_{\chance,j}
    \end{array}
    \right.\hspace{-0.5em}
    \right\},
\end{equation}
where $b_{\chance,j}$, $E_{\chance,j}$ and $F_{\chance,j}$
have appropriate dimensions, ${r_j\in\N_0}$, $\mathcal{K}_{\chance,j}$ is a cone
and it is implied that $\mathcal{A}_{\chance, j}$ is a subset of the probability
simplex of $\R^{n_j}$  \cite{SSP2019a}.

For example, for $\AVaR_{\chance}$ it is
\(
{E_{\chance, j} = [I_{n_j}~-I_{n_j}~1_{n_j}]\tran}
\),
$r_j=0$,
\(
b_{\chance, j}=[\chance^{-1}\pi_j\tran~0~1]\tran
\)
and
\(
\mathcal{K}_{\chance, j} = \R^{2n_j} \times \{0\}
\).
From convex duality, if there exist \(\mu^*\in\R^{n_j}\) and
\(\nu^*{}\in{}\R^{r_j}\) so that \(b - E\mu^* - F\nu^*\) is in the relative
interior of $\mathcal{K}_j$, which is a very weak assumption,
then $\rho_j[C_j]$, for $C_j\in\R^{n_j}$ can be written as
\begin{equation}
    \label{eq:risk-dual-representation}
    \rho_j[C_j]
    {}={}
    \min_{y}
    \left\{
    y\tran b_{\chance, j}
    \left|
    \hspace{-0.3em}
    \begin{array}{l}
        E_{\chance, j}\tran y \,{=}\, C_j,
        \,
        F_{\chance, j}\tran y \,{=}\, 0,
        \\
        y {}\in{} \mathcal{K}_{\chance, j}^*
    \end{array}
    \right.\hspace{-0.6em}
    \right\},
\end{equation}
where \(\mathcal{K}_j^*\) is the convex dual of \(\mathcal{K}_j\) \cite[Thm~2.4.1]{BTN2001}.
By virtue of \eqref{eq:risk-dual-representation}, constraints
\eqref{eq:stagewise-risk-constraints} are equivalent to the existence
of a
\(
y_j \in \mathcal{K}_{\chance, j}^*
\)
such that
\(
E_{\chance, j}\tran y_j = C_j
\)
and
\(
F_{\chance, j}\tran y_j = 0.
\)
This makes Problem \eqref{eq:optimal-control-problem} into a mixed-integer
conic optimization problem, yet without additional binary variables as it
was the case with chance constraints.


\section{Simulations}

\begin{figure}
  \centering
  \includegraphics{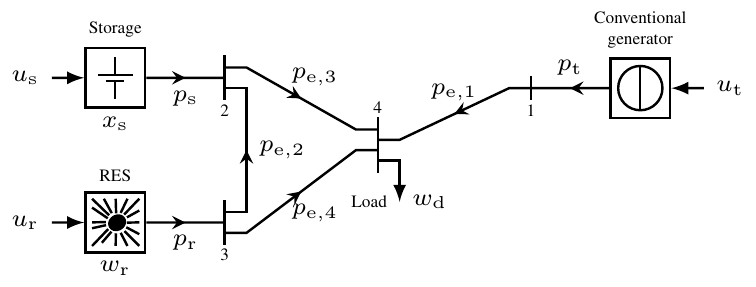}%
  \caption{Simple microgrid with storage unit, renewable and conventional generator as well as load. Source: \cite{HSR+2020}.}
  \label{fig:testSetupWithLoad}
\end{figure}

In the closed-loop simulations, the grid in \cref{fig:testSetupWithLoad} was considered.
It is composed of a \ac{pv} power plant, a conventional generator, and a storage unit, i.e., all basic components from \cref{sec:microgrid-model} are included.
The units are connected to each other and to a load via transmission lines that all have a rated power of \unit[1.3]{pu}.
For all transmission lines a susceptance of \unit[-20]{pu} was considered.
The remaining parameters of the units and the weights of the cost function can be found in \cref{tab:modelSpecificVariables}.

\begin{table}
  \centering
  \newcommand{\cIn}[2]{\ensuremath{c_{\mathrm{#1}}^{#2}}}
  \newcommand{\cS}[1]{\cIn{s}{#1}}
  \newcommand{\cT}[1]{\cIn{t}{#1}}
  \newcommand{\cR}[1]{\cIn{r}{#1}}
  \caption{Unit parameters and weights of cost function.}
  \label{tab:caseStudy:ModelParameters}
  \begin{tabular}{clc@{}cl}
    \toprule
    Parameter                                                                     & Value                     &  & Weight              & Value                               \\
    \cmidrule{1-2} \cmidrule{4-5}
    $[{p}_{\mathrm{t}}^{\min}, {p}_{\mathrm{r}}^{\min}, {p}_{\mathrm{s}}^{\min}]$ & $[0.4, 0, -1]\,\unit{pu}$ &  & $\cT{}$             & $\unit[0.1178]{}$                   \\[1pt]
    $[{p}_{\mathrm{t}}^{\max}, {p}_{\mathrm{r}}^{\max}, {p}_{\mathrm{s}}^{\max}]$ & $[1, 2, 1]\,\unit{pu}$    &  & $\cT{\prime}$       & $\unit[0.751]{\nicefrac{1}{pu}}$    \\[1pt]
    $[{x}^{\min}, {x}^{\max}]$                                                    & $[0, 4]\,\unit{pu\,h}$    &  & $\cT{\prime\prime}$ & $\unit[0.0693]{\nicefrac{1}{pu^2}}$ \\[1pt]
    $[\tilde{x}^{\min}, \tilde{x}^{\max}]$                                        & $[1, 3]\,\unit{pu\,h}$    &  & $\cT{\mathrm{s}}$   & $\unit[0.3162]{}$                   \\[1pt]
    $x^0$                                                                         & $\unit[3]{pu\,h}$         &  & $\cR{}$             & $\unit[1]{\nicefrac{1}{pu}}$        \\[1pt]
    $[{K}_\mathrm{t}, {K}_\mathrm{s}]$                                            & $[1, 1]$                  &  & $\gamma$            & $\unit[0.95]{}$                     \\
    \bottomrule
  \end{tabular}
\end{table}

The controllers and the \ac{mg} model in \cref{fig:testSetupWithLoad} were implemented in MATLAB~R2019b using YALMIP~R20200116~\cite{Lof2004} and solved with Gurobi~9.0.2 on a machine with a \unit[3.70]{GHz} Intel$^\text{\textregistered}$ Xeon$^\text{\textregistered}$ E5-1620 v2 CPU and \unit[32]{GB}\,RAM.
The computation times were reduced by using the results of the previous iterations for a warm-start of the solver.
Moreover the binary variables used in the model were relaxed for stages larger than or equal to $j=4$ to speed up the solver.
Note that the binary variables used to formulate the chance constraints in \eqref{eq:chance-constrained-optimal-control-problem} were not relaxed.

The time series of load demand used in the simulations was based on measurements from a real-world islanded \ac{mg}.
For the \ac{pv} power plant, irradiance data from \cite{ARM2009} was employed.
The time series of load and irradiance were used to train seasonal \ac{arima} models \cite{BJR2013} that are employed to forecast load and available renewable infeed for the \ac{mpc} formulation.
In detail,
for load an $\operatorname{ARIMA}(10, 0, 8)(7, 1, 7)_{48}$ model,
and for irradiance an $\operatorname{ARIMA}(6, 1, 2)(1, 1, 1)_{48}$ model were used.
From both models, scenario trees of the prediction errors were constructed in a similar way as in \cite{HR2003,HSB+2015}.
In what follows, chance constraints and the corresponding risk constraints are imposed
separately for the upper and lower bounds, that is, we use
\(
d_j^{\Max}(x)
{}={}
\min_{y \leq \tilde{x}_{s}^{\Max}}|x-y|,
\)
and
\(
d_j^{\Min}(x)
{}={}
\min_{y \geq \tilde{x}_{s}^{\Min}}|x-y|.
\)

\subsection{Closed-loop simulations for $\alpha = 0.5$}

\begin{figure}
  \centering
  \includegraphics{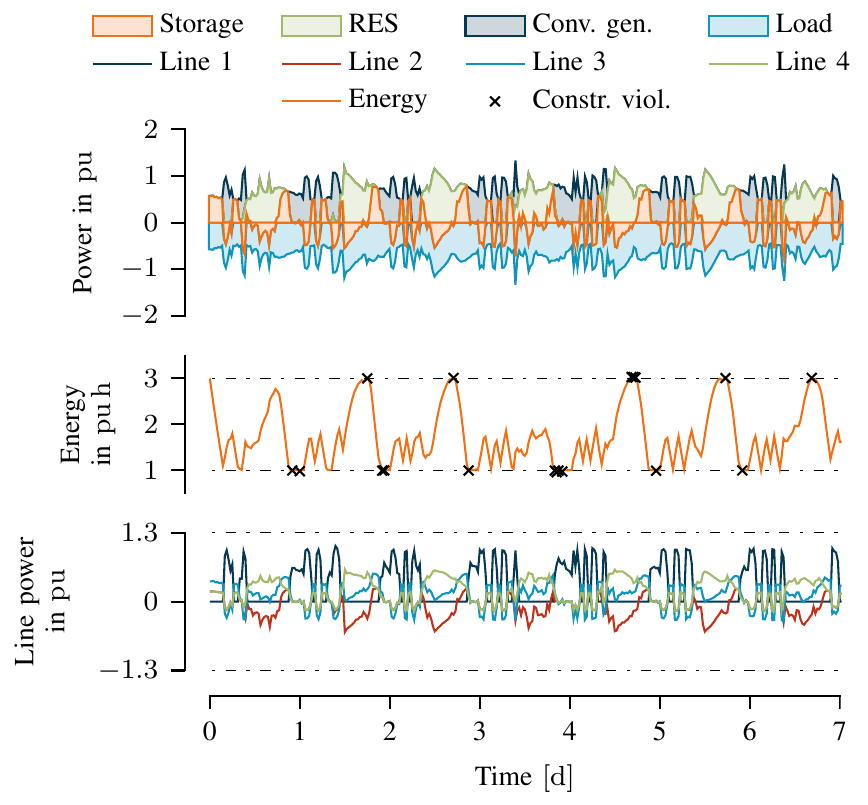}
  \caption{Results of closed-loop simulations over one week with risk-constrained \ac{mpc} at $\alpha = 0.5$.}
  \label{fig:caseStudyResults}
\end{figure}

The results of the closed-loop simulations with the risk-constrained \ac{mpc} ($\alpha = 0.5$) over a simulation horizon of \unit[7]{d}, i.e., 336 simulation steps, are shown in \cref{fig:caseStudyResults}.
It can be noted that during each day, \ac{pv} infeed is used to charge the storage unit.
During this time, the conventional unit is disabled.
At the end of each day, the storage unit is discharged as the infeed from the \ac{pv} power plant decreases.
At night, the conventional unit is repeatedly enabled to provide power to the loads and charge the storage units.
In theory, $\alpha = 0.5$ allows for \unit[50]{\%} of the predicted values to be above $\tilde{x}^{\max}=\unit[3]{pu\,h}$ at each stage and for \unit[50]{\%} of the values to be below $\tilde{x}^{\min} = \unit[1]{pu\,h}$.
In the closed-loop simulations, however, only $18$ energy values outside the interval $[\tilde{x}^{\Min}, \tilde{x}^{\Max}]$ were observed with maximum distance to the interval below \unit[0.05]{pu\,h}.

\subsection{Comparison of closed-loop simulations for different $\alpha$}

\begin{table}[t]
  \centering
  \caption{Accumulated values of closed-loop simulation with simulation horizon $K = 336$.}
  \label{tab:CaseStudyComparison}
  \newcommand{\centerMathNum}[1]{\noDecimal{\parbox{5mm}{\centering \ensuremath{#1}}}}
  \newcommand{\csvDispColumn}[1]{%
    \csvreader[head to column names, separator=comma, filter equal={\thecsvinputline}{3}]{data/simulationResults.csv}{}{#1} &
    \csvreader[head to column names, separator=comma, filter equal={\thecsvinputline}{4}]{data/simulationResults.csv}{}{#1} &
    \csvreader[head to column names, separator=comma, filter equal={\thecsvinputline}{5}]{data/simulationResults.csv}{}{#1} & &
    \csvreader[head to column names, separator=comma, filter equal={\thecsvinputline}{6}]{data/simulationResults.csv}{}{#1} &
    \csvreader[head to column names, separator=comma, filter equal={\thecsvinputline}{7}]{data/simulationResults.csv}{}{#1} &
    \csvreader[head to column names, separator=comma, filter equal={\thecsvinputline}{8}]{data/simulationResults.csv}{}{#1}
  }
  \begin{tabular}{rd{0.0}d{0.0}d{0.0}c@{}d{0.0}d{0.0}d{0.0}}
    \toprule
                                                    & \multicolumn{3}{c}{Chance-constr., $\alpha=$} &                     & \multicolumn{3}{c}{Risk-constr., $\alpha=$}                                                                      \\
    \cmidrule{2-4}  \cmidrule{6-8}
                                                    & \centerMathNum{0.1}                           & \centerMathNum{0.2} & \centerMathNum{0.5}                         &  & \centerMathNum{0.1} & \centerMathNum{0.2} & \centerMathNum{0.5} \\
    \midrule
    Avg. costs $\bar{\ell}$                         & \csvDispColumn{\costs}                                                                                                                                                                 \\
    Avg. ren. share $[\unit{\%}]$                   & \csvDispColumn{\renewableShare}                                                                                                                                                        \\
    \midrule
    $x \notin [\tilde{x}^{\min}, \tilde{x}^{\max}]$ & \csvDispColumn{\errEnergy}                                                                                                                                                             \\
    Switching actions                               & \csvDispColumn{\switchingActions}                                                                                                                                                      \\
    \midrule
    Avg. solve time $[\unit{s}]$                    & \csvDispColumn{\meanSolverTime}                                                                                                                                                        \\
    Max. solve time $[\unit{s}]$                    & \csvDispColumn{\maxSolverTime}                                                                                                                                                         \\
    \bottomrule
  \end{tabular}%
\end{table}%

The results of closed-loop simulations with the chance-constraint \ac{mpc} \eqref{eq:chance-constrained-optimal-control-problem} and the risk-constrained \ac{mpc} \eqref{eq:risk-constrained-optimal-control-problem} for different values of $\alpha$ are shown in \cref{tab:CaseStudyComparison}.
It can be noted that with increasing $\alpha$, the average closed-loop costs
\[
  \bar{\ell} = \textstyle \sum\limits_{k=1}^{336} \ell(v(k), v(k-1), p(k)),
\]
decrease.
Comparing, for example, the risk-constrained \ac{mpc} for $\alpha = 0.1$\footnote{Note that a uniform imposition of constraints via the chance-constrained \ac{mpc} with $\alpha = 0$ also leads to an average cost of $\bar\ell = 3.36$.} and $\alpha = 0.5$, one can see that the average cost decreases about \unit[1]{\%}.
Note that this decrease is reached solely by allowing some energy constraints violations.
As noted earlier, the maximum violation of the chance constraints is less than \unit[0.05]{pu\,h},
so the decrease in price comes at the acceptable disadvantage of very small violations of $x$.

It can be seen in \cref{tab:CaseStudyComparison} that the number of constraint violations of the chance-constrained approach is smaller than that of the risk-constraint approach.
Moreover, the average closed-loop cost $\bar{\ell}$ of the chance-constrained \ac{mpc} is slightly smaller than that of the risk-constrained approach.
Both effects are based on the fact that risk constraints overapproximate chance constraints.
Furthermore, the average and the maximum computing time of the solver are lower for the risk-constrained approach.
Unlike chance constraints, the risk-based aporoach does not neccesitate additional binary varibles.
This leads to optimization problems that can accomodate an adequate number of storage units at a reasonable computational cost.

\begin{figure}
  \centering
  \includegraphics{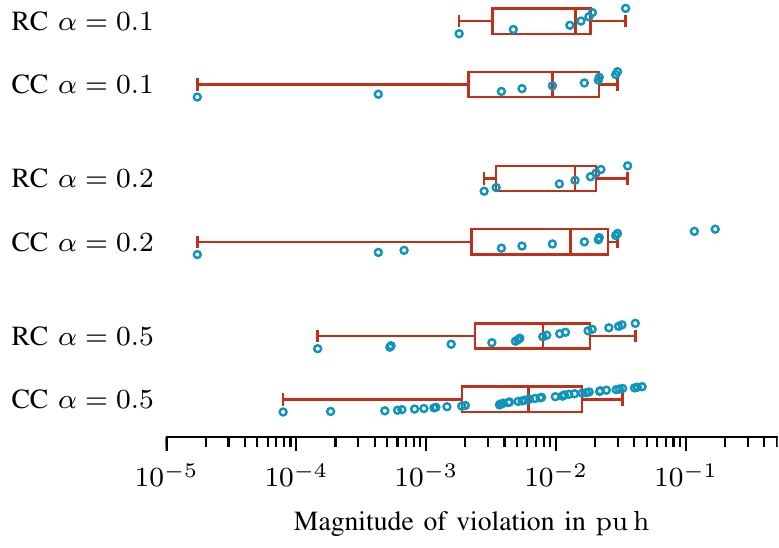}
  \caption{Distribution of constraint violations relative to bounds, i.e., ${\tilde{x}^{\min} - x(k)}$ and $x(k) - \tilde{x}^{\max}$ for chance constrained (CC) and risk constrained (RC) \ac{mpc}. Note that states that did not violate any constraint were omitted.}
  \label{fig:boxplotConstraintViolations}
\end{figure}

In \cref{fig:boxplotConstraintViolations}, distributions of values outside the interval ${[\tilde{x}^{\min}, \tilde{x}^{\max}]}$ are shown for the different \ac{mpc} approaches and  values of $\alpha$.
It can be noted that, for the same value of $\alpha$, the risk-constrained controller typically comes with a smaller number of violations.
The largest violations of all simulations occur using the chance-constrained controller.
This indicates that imposing risk constraints can help to prevent extreme violations by taking their magnitude into account.


\section{Conclusions}
This paper proposes a scenario-based stochastic model predictive 
control formulation with stagewise risk-based constraints on the 
state of charge of the storage units of a microgrid. 
Unlike chance constraints, risk constraints can control both the 
frequency and magnitude of constraint violations; 
they can allow infrequent and mild violations of the state-of-charge
constraints leading to a lower operating cost.
Lastly, risk constraints are convex unlike chance constraints, which 
necessitate the introduction of binary variables; as a result, 
risk-constrained problems can be solved faster.

\bibliographystyle{IEEEtran}
\bibliography{literature}
\end{document}